\documentclass[12pt,a4paper]{article}
\usepackage[english,french]{babel}
\usepackage[latin1]{inputenc}

\usepackage[centertags]{amsmath}
\usepackage{amsfonts}
\usepackage{amssymb}
\usepackage{amsthm}
\usepackage{newlfont}
\usepackage{color}      
\usepackage{epsfig}
\usepackage{graphicx}
\usepackage{float}

\newtheorem{theo}{Theorem}
\newtheorem{prop}{Proposition}
\newtheorem{defi}{Definition}

\newtheorem{lem}{Lemma}

\newenvironment{rem}{\noindent \textbf{Remark.}}{\medskip}

\newcommand{\complex}{\mathbb C}
\newcommand{\integer}{\mathbb N}
\newcommand{\rinteger}{\mathbb Z}

\newcommand{\I}{\mathcal{I}}
\newcommand{\J}{J}
\newcommand{\Y}{\mathfrak{Y}}
\newcommand{\val}{\nu}
\renewcommand{\tild}[1]{\widetilde{#1}}

\newcommand{\pmatrice}[4]{\begin{pmatrix}#1&#2\\#3&#4 \end{pmatrix}}

\title{Dynamics of rational symplectic
mappings and difference Galois theory}

\author{
Guy Casale\thanks{ {\bf address} : IRMAR UMR 6625, Universit\'e de Rennes 1, Campus de Beaulieu 35042 Rennes Cedex - France ;\quad 
{\bf
E-mail} : {\sf guy.casale@univ-rennes1.fr }} 
 \ $\&$ Julien Roques\thanks{ {\bf address} :\'Ecole Normale
Sup\'erieure,  D\'epartement de Math\'ematiques et
  Applications UMR 8553, 45, rue d'Ulm, 75230 Paris Cedex 05 - France 
 ;\quad 
 {\bf E-mail} :
\textsf{julien.roques@ens.fr} } 
}

\begin{document}
\selectlanguage{english}

\maketitle

\hrule
\vspace{10pt}

\begin{abstract}
\textit{ In this paper we study the relationship between the
integrability of rational symplectic maps and difference Galois
theory. We present a Galoisian condition, of Morales-Ramis type,
ensuring the non-integrability of a rational symplectic map in the
non-commutative sense (Mishchenko-Fomenko). As a particular case, we obtain a complete discrete analogue of Morales-Ramis Theorems for non-int\-egra\-bi\-li\-ty in the sense of Liouville.}
\end{abstract}

\selectlanguage{francais}
\begin{abstract}
{ \begin{center} \textsc{Dynamique des applications rationnelles
symplectiques\\ et th\'eorie de Galois aux diff\'erences \\ \medskip }\end{center}

\textit{Dans cet article nous \'etudions les relations entre
l'int\'egrabilit\'e des applications rationnelles symplectiques et la
th\'eorie de Galois des \'equations aux diff\'erences. Nous donnons une
condition galoisienne du type Morales-Ramis de
non-int\'egrabilit\'e d'une application rationnelle
symplectique au sens non-commutatif de Mishchenko-Fomenko
 et en particulier au sens de Liouville .} }
\end{abstract}

\selectlanguage{english}

\vspace{10pt}
 \hrule
 \vspace{20pt}

{\small
  The first author is supported by ANR project ``{\it Int\'egrabilit\'e en
m\'ecanique hamiltonienne}'' n$^\circ$JC$05\_41465$.
  The second one is supported by the ANR project 
``{\it Ph\'enom\`ene de Stokes, renormalisation, th\'eories de Galois}''. \\

}

 \newpage

\tableofcontents
\normalsize

\medskip

\newpage

\section{Introduction-Organization}
\label{introduction-organisation}

The problem of deciding whether a given continuous dynamical
system is integrable is an old and difficult problem. Many
examples arise in classical mechanics and physics, the \textit{three body problem} being one of the most famous. In 1999, after anterior
works by a number of authors among whom Kowalevskaya, Poincar\'e,
Painlev\'e and more recently Ziglin \cite{ziglin}, Morales and Ramis exhibited
in \cite{morales,moralesramis1,moralesramis2} an algebraic obstruction to the
integrability of Hamiltonian dynamical systems. This obstruction
lies in the differential Galois group of the linearization of the
dynamical system under consideration along a given particular
solution (the variational equation). During the last years, the
method developed by Morales and Ramis have shown its efficiency,
permitting the resolution of numerous
classical problems. \\

Similar questions occur in the study of discrete dynamical
systems. For instance, a rational map $f:\mathbb{P}^2(\complex)
\rightarrow \mathbb{P}^2(\complex)$ being given the typical
question is : does there exist a non-trivial rational invariant of
$f$ ? That is, does there exist $H\in\complex(x,y)$ non-constant
such that $H\circ f = H$ ? In case of a positive answer, $f$ is
called integrable (we refer to section \ref{integrabilite} for the
general notion of integrability). Inspired by Morales-Ramis
theory, we point out an algebraic criterion of non-integrability
based on difference Galois theory. Our method involves the
linearization of the dynamical system under consideration along a
``particular solution" (that is an invariant parametrized complex
algebraic curve of genus zero
$\iota\left(\mathbb{P}^1(\complex)\right)$) called the discrete
variational equation. From a technical point of view, we present
two proofs of our main theorem : the first one is based on the
difference Galois theory developed by Singer and van der Put in
\cite{singer-vanderput} (there are some complications with respect
to the differential case due to the fact that the Picard-Vessiot
ring associated to a given difference equation is not necessarily
a domain), whereas the second one is based on
Malgrange theory.\\

The paper is organized as follows. In sections \ref{section
Notations Hypotheses} and \ref{integrabilite}, definitions and
notations used in the article are given. In section
\ref{nonlineaire} we introduce the Galois groupoid of B. Malgrange
and prove that for an integrable system the groupoid of Malgrange
is infinitesimally commutative. But commutativity of the Malgrange
groupoid is not easily checkable.
 In section \ref{lineaire}, we recall some facts about difference equations from \cite{singer-vanderput}. Section \ref{linearisation} deals
with le linearization of the discrete system and its differential
invariants along an invariant curve. The main theorem is proved in
section \ref{theoreme} :

\begin{theo}
Let $f$ be an integrable symplectic map having a
$\phi$-adapted curve $\iota$. Then the neutral component of the
Galois group of the discrete variational equation of $f$ along
$\iota$ is commutative.
\end{theo}

The two proofs are given in this section. The first one deals
directly with the linearized system following the original
approach of Morales-Ramis \cite{moralesramis1}. The second used
the linearization to check the commutativity of Malgrange
groupoid. This work is based on a suggestion of Ramis; we are pleased to thank him.

\section{Notations-Hypotheses}\label{section Notations Hypotheses}
Unless explicitly stated otherwise, in the whole article we will
use the following notations and conventions :
\begin{itemize}
\item[$\bullet$] the (symplectic) manifolds and mappings are complex;
\item[$\bullet$] $V$ will denote a complex algebraic variety of even complex dimension $2n$ ($n\in \integer^*$);
\item[$\bullet$] $\complex(V)$ will denote the field of rational functions over the algebraic variety $V$;
\item[$\bullet$] $\omega$ will denote a non-degenerate closed holomorphic $2$-form on $V$ (symplectic form);
\item[$\bullet$] $f : V \dasharrow V$ will denote a
birational mapping preserving $\omega$ that is such that $f^*\omega
= \omega$ (symplectic mapping of $(V,\omega)$);
\item[$\bullet$] the symplectic gradient of a rational function $H : V \dasharrow \complex$ is the rational vector field $X_{H}$ satisfying
$dH = \omega(X_{H} ,\, \cdot \, )$;
\item[$\bullet$] the induced Poisson bracket on $\complex(V)$ will be denoted by  $\{ H_{1} ,H_{2} \} = \omega(X_{H_{1}}, X_{H_{2}})$.
\end{itemize}

\section{Integrability of a symplectic map}
\label{integrabilite}

 For the sake of completeness, we first
recall the notion of Liouville integrability
of a continuous hamiltonian dynamical system.

\begin{defi}[\cite{liouville}]

Let $H$ be a rational Hamiltonian on $V$. It is said to be integrable in the Liouville sense (or in the Langrangian way) if there is $H_{1}, \dots , H_{n}$ functionally independent first integrals of $H$ (\textit{i.e.} $\{H , H_{j}\} = 0$) in involution :

\begin{itemize}
\item  functionally independent means that their differentials $dH_1(z)$, ..., $dH_n(z)$ are linearly
independent for at least one point $z \in V$ or, equivalently, for
all $z$ in a Zariski dense open subset of $V$;
\item[$\bullet$] in involution means that, for all
  $i,j$ in $[[1,n]]$,
  $\{H_i,H_j\}=0$.
\end{itemize}
\end{defi}

This means that the motion takes place in the fibers of a Lagrangian fibration. This notion of integrability was generalized by Fomenko-Mishchenko \cite{fomenko-mishchenko} (see also \cite{nekhoroshev, bogoyavlenskij})
following previous work of Marsden-Weinstein \cite{MW}.
\begin{defi}
Let $H$ be a rational Hamiltonian on $V$. It is said to be integrable in the non-commutative sense (or in the isotropic way) if there is $H_{1}, \dots , H_{n + \ell}$ independent first integrals of $H$ \textit{i.e.} $\{H , H_{j}\} = 0$ such that the distribution $\bigcap_{i} \ker dH_{i}$ is generated by the symplectic gradients $X_{H_{1}}, \dots, X_{H_{n-\ell}}$ as vector space over $\complex(V)$.
\end{defi}

This means that the motion takes place in the fibers of an isotropic fibration. Mishchenko and Fomenko conjectured that this notion
of integrability is equivalent to Liouville's one. This question will not be discuss here. For a proof in the finite dimensional case, we refer
the reader to \cite{sadetov, bolsinov}. A Galoisian condition for integrability in the isotropic way was obtained by A. Maciejewki and
M. Przybylska in \cite{MP}. As an evidence of Mishchenko-Fomenko conjecture, the Galoisian conditions in both cases are the same.

In the context of discrete dynamical systems (iteration of maps)
the same notion of integrability is used.

\begin{defi}
The symplectic map $f$ is integrable (in the isotropic way)
if it admits $n+\ell$ first integrals $H_1,...,H_{n+ \ell}$ in
$\complex(V)$, functionally independent such that
$\bigcap_{i} \ker dH_{i}$ is generated by the symplectic gradients $X_{H_{1}}, \dots, X_{H_{n-\ell}}$

\end{defi}

The only point in this definition requiring some comments is the
notion of first integral. A rational function $H \in \complex(V)$ is a first integral of $f$
if $H$ is constant over the trajectories of $f$, that is $H \circ
f =H$.

\section{A non-linear approach}
\label{nonlineaire}

\subsection{Malgrange groupoid and its Lie algebra}

Let $V$ be a dimension $n$ smooth complex algebraic variety. Let $R_{q}V$ denote the order $q$ frame bundle on $V$. This bundle is the space of order $q$ jets of invertible maps  $r : (\complex^n,0) \to V $. Such a jet will be denoted by $j_{q}(r)$.

As soon as one gets local coordinates $r_{1}, \ldots , r_{n}$ on an open subset $U$ of $V$, one gets local coordinates $r_{i}^\alpha$ with $i \in [[1, \ldots, n]]$ and $\alpha \in \integer^n$ such that  $|\alpha| \leq q$ on the open set of $R_{q}V$ of order $q$ frames with target in $U$.
This bundle is an algebraic variety and is naturally endowed with two algebraic actions. The first one is the action of the algebraic group $\Gamma_{q}^n$ of order $q$ jets of invertible maps $(\complex^n,0) \to (\complex^n,0)$ by composition at the source. For this action $R_{q}V$ is a $\Gamma_{q}^n$-principal bundle on $V$.

The second action is the action of the algebraic {\it groupoid} $J^*_{q}(V,V)$ of invertible local holomorphic maps $s : (V,a) \to (V,b)$ by composition at the target.

Let us recall the definition of an algebraic groupoid on $V$.
\begin{defi}
An algebraic groupoid on $V$ is an algebraic variety $G$ with :
\begin{itemize}
\item[$\bullet$] two maps $s,t : G \to V$, the {\it source} and the {\it target},
\item[$\bullet$] a associative composition $\circ : G \underset{V}{_s \times_t} G \to G$,
\item[$\bullet$] an identity $id: V \to G$ such that for all $g \in G$
$$ g \circ id(s(g)) = id(t(g)) \circ g = g,$$
\item[$\bullet$] an inverse ${}^{\circ-1} : V \to V$ such that $\forall g \in G$
$$ \  s(g^{\circ-1}) = t(g) \ ,  \ t(g^{\circ-1}) = s(g) \ , \  g \circ g^{\circ-1} = id(t(g)) \ \text{and} \   g^{\circ-1} \circ g = id(s(g)).$$
\end{itemize}
The induced morphisms on the structure (sheaves of) rings must satisfy the dual diagrams.
\end{defi}

The varieties $J^*_{q}(V,V)$ are archetypes of groupoid on $V$.

\begin{defi}
Let $f :  V \dasharrow V$ be a rational map. The order $q$ prolongation of $f$ is $R_{q}f : R_qV \dasharrow R_qV$ defined by $(R_{q}f) (j_{q}(r)) = j_{q}(f \circ r)$.
\end{defi}
The basic objects used to understand the dynamic of $f$ are (rational) first integrals also called (rational) invariants, {\it i.e.} functions $H : V \dasharrow \complex$ such that $H\circ f = H$. From a differential-algebraic point of view, invariants of $f$ are replaced by differential invariants of $f$.
\begin{defi}
A differential invariant of order $q$ of $f$ is an invariant of $R_{q}f$, {\it i.e.} $H : R_{q}V \dasharrow \complex$ such that $H\circ R_{q}f = H$.
\end{defi}
Thanks to the canonical projections $\pi_{q}^{q+1} : R_{q+1}V \to R_{q}V$, an order $q$ invariant is also an order $q+1$ invariant. Let us denote by $\mathcal{D}_{q}Inv(f)$ the set of all the differential invariants of order $q$ of $f$. This is a $\complex$-algebra of finite type.
Invariant tensor fields are particular cases of order $1$ differential invariants.

\begin{defi}
Let $V$ be a smooth complex algebraic variety and  $f :  V \dasharrow V$ be a rational map.
The  $\mathcal{D}_{q}$-closure of $f$, denoted by $G_{q}(f)$, is the subvariety of $J^*_{q}(V,V)$ defined by the equations
$H\circ j_{q}(s) = H $ for all $H \in \mathcal{D}_{q}Inv(f)$.
\end{defi}

\begin{rem}
In general,  these subvarieties are not subgroupoids but subgroupoids above a Zariski dense open subset of $V$.
\end{rem}

Let $J^*(V,V)$ be the projective limit $\underset{\leftarrow}{\lim}\,J_{q}(V,V)$ and $\pi_{q}$ be the projection of this space on $J_{q}(V,V)$.

\begin{defi}[Malgrange]
Let $V$ be a smooth complex algebraic variety and  $f :  V \dasharrow V$ be a rational map.
The Malgrange groupoid of $f$ is the subvariety of $J^*(V,V)$ given by $\bigcap \pi_{q}^{-1} G_{q}(\Psi)$ and it is denoted by $G(f)$.
\end{defi}

\begin{rem}
This is not the original definition from Malgrange
\cite{malgrange} which is about foliation. The extension of the
definition to discrete dynamical systems is straightforward.
Proposition 2.36 from chapter 3 of \cite{pommaret} (pp 467--469)
allows us to `simplify' the exposition by using directly
differential invariants.
\end{rem}

A solution of $\bigcap \pi_{q}^{-1} G_{q}(\Psi)$ is a formal map $s : \widehat{V,a} \to \widehat{V,b}$ whose coefficients satisfy all the equations of the $G_{q}(\Psi)$. It is said to be convergent if the formal map converges. There are well-known solutions of $G(f)$:  the iterates $f^{\circ n}$ preserve all the differential invariants of $f$.

For a rational map $\Phi : \mathbb{P}^1(\complex) \to \mathbb{P}^1(\complex)$, the Malgrange groupoid is computed in \cite{casale-toulouse}. Following this computation the maps with a `small' Malgrange groupoid are the maps called integrable {\it i.e.} maps with a non trivial commutant \cite{ritt, veselov,buiumzimmerman}.

Groupoid are not amenable objets. It is much more easy to work with the infinitesimal part of a groupoid : its (analytic sheaf of) Lie algebra. For subgroupoids of $J^*(V,V)$ these Lie algebras are algebras of vector fields.
Let us give some notations before the definitions.
Let $r_{1}, \ldots r_{n}$ be local coordinates on $V$ and $r_{i}^\alpha$ for $1\leq i \leq n$ and $\alpha \in \integer^n$ with $|\alpha| \leq q$ be the induced coordinates. A function on $R_{q}V$ can be differentiated with respect to coordinates on $\complex^n$ by mean of the chain rule {\it i.e.} using the vectors fields
$$
D_{j} =\sum_{i,\alpha} r_{i}^{\alpha + \epsilon(j)}\frac{\partial}{\partial r_{i}^{\alpha}}
$$
where $\epsilon(j) = \overset{\stackrel{j}{\downarrow}}{(0\ldots 0,1,0\ldots 0)}$, one gets a function on $R_{q+1}V$.

\begin{defi}
Let $X = \sum_{i}  a_{i} \frac{\partial}{\partial r_{i}}$ be a holomorphic vector field on an open subset $U \subset V$. The order $q$ prolongation of $X$ is the vector field $R_{q}X = \sum_{i,\alpha}  D^\alpha a_{i}  \frac{\partial}{\partial r^{\alpha}_{i}}$ on $R_{q}U \subset R_{q}V$.
 \end{defi}

 \begin{defi}
Let $V$ be a smooth complex algebraic variety and  $f :  V \dasharrow V$ be a rational map.
The Lie algebra of the Malgrange groupoid of $f$ is the analytic sheaf of vector field $X$ on $V$ such that
$$
\forall q \in \integer \text{ \ and \ } \forall H \in \mathcal{D}_{q}Inv(f), \ (R_{q}X)H = 0.
$$
 \end{defi}

 \begin{rem}
The original definition from \cite{malgrange} is richer about the structure of this sheaf but for the sequel this definition will be sufficient.
 \end{rem}

\subsection{Malgrange groupoid of an integrable symplectic map}
\label{invariantsdifferentiels}
As we already said invariant tensors are examples of order $1$ differential invariants.
We will describe these invariants and their lifting up in some cases.

\begin{lem}
Let $T$ be a rational tensor field on $R_q V$. It induces $m$
functions $h_T^{i}$, $ 1\leq i \leq m$, on $R_{q+1}V$. Furthermore
for a vector field $X$ on $V$ one gets
$$
L_{R_q X} T  = 0 \text{ if and only if } (R_{q+1}X) h_T^{i} = 0 \quad \forall 1\leq i \leq m.
$$
\end{lem}

Let's see on examples these functions and this property.\\

\noindent
{\bf Rational function --}  A rational function $h : V \dasharrow \complex$ can be lifted to a rational function on $R_{q}V$. It is clear that if $X$ is a vector field on $V$ and $h$ is $X$-invariant then its lifting up is $R_{q}X$-invariant. \\

\noindent
{\bf Rational form --}  By the very definition, if $\omega$ is a $p$-form on $V$ then it is a function on $(TV)^p = J_{1}((\complex^p,0) \to V)$. So it induces $n \choose p$ functions on $R_{1}V$. One can then lift these functions on higher order frame bundles.

Now if $X$ is a vector field on $V$ and $\omega$ is a $X$-invariant $p$-form. Let $r_{1}, \ldots, r_{n}$ be local coordinates on $V$,
one gets
$$
\omega = \sum_{I \in C(p,n)} w_{I}dx_{\wedge I} \text{ and } X = \sum_{i=1}^n a_{i} \frac{\partial}{\partial r_{i}}
$$
 where $C(p,n)$ stands for the set of length $p$ increasing sequences
$i_1 < \ldots < i_p$ in $\{1,\ldots, n\}$ and $dr_{\wedge I}$ stands for $dr_{i_1}\wedge \ldots \wedge dr_{i_p}$.
The lifting up of $\omega$ to $R_{1}V$ is the $n \choose p$-uple of functions :
$$
h_\omega^{J}(- r_{i} - r_{i}^{\epsilon(j)} - ) = \sum w_{I} r^{J}_{I}
$$ where $J \subset \{1,\ldots, n\}$ with $\# J = p$ and $r^{J}_{I}$ stands for the determinant of the matrix $\begin{pmatrix}
r_{i}^{\epsilon(j)}\end{pmatrix}_{i \in I \atop j \in J}$.

Easy computations give
$$
\begin{array}{rcl}
\displaystyle L_X \omega & = & \displaystyle \sum X w_I dr{\wedge I} + \sum (-1)^{\ell + 1} w_I \frac{\partial a_p}{\partial r_k} dr_k \wedge dr_{\wedge I-\{i_\ell\}}, \\
\displaystyle R_1 X & = & \displaystyle \sum a_i \frac{\partial}{\partial r_i} + \sum \frac{\partial a_i}{\partial r_k} r_k^{\epsilon(j)}\frac{\partial}{\partial r_i^{\epsilon(j)}},\\
\text{ and }\displaystyle  (R_1X)(h_\omega ^ J) & = &\displaystyle  \sum Xw_I r_I^J + \sum w_I \frac{\partial a_i}{\partial r_k} r_k^{\epsilon(j)} \frac{\partial r_I^J}{\partial r_i^{\epsilon(j)}}.\\
\text{ Thus one has }\displaystyle  h_{L_X \omega}^J & = & \displaystyle (R_1X)(h_\omega^J).
\end{array}
$$
For instance, starting from a $X$-invariant function $h$ on $V$ one gets $n + 1$ order $1$ differential invariants : $h$ itself and $h_{dh}^1, \ldots , h_{dh}^n$.\\

\noindent
{\bf Rational vector field --} Let $v$ be a rational vector field on $V$ : $v = \sum v_i \frac{\partial}{\partial r_i}$ in local coordinates. An order $1$ frame
determines a basis $e_1, \ldots, e_n$ of $T_rV$. The coordinates of $v(r)$ in this basis determine $n$ functions $h_v^i : R_1 V \dasharrow \complex$. \\

\noindent
{\bf Higher order forms --} A $p$-form $\omega$ on $R_q V$ determines functions on $R_1 R_q V$. The induced function on $R_{q+1}V$ are the
 pull-back by the natural inclusion $R_{q+1}V \subset R_1 R_q V$.

\begin{theo}
\label{malgcommutatif}
The Lie algebra of Malgrange groupoid of an integrable (in the isotropic way) symplectic map is commutative.
\end{theo}
\begin{proof}
Vector fields of the Lie algebra of Malgrange groupoid  must
preserve $n + \ell$ independent functions $H_{1}, \dots , H_{n +
\ell}$:  $X H_{i} = 0$. Because they preserve $H_{i}$ and the
symplectic form $\omega$, they must preserve the symplectic
gradients $X_{{i}}$.

On the subvarieties  $H_{i}=c_{i}$ the integrable map preserves a parallelism given by $X_{H_{1}}, \dots, X_{H_{n-\ell}}$.
This parallelism is commutative : $[X_{H_{i}}, X_{H_{j}}] = X_{H_{j}} H_i= 0$ for all $1 \leq i \leq n+ \ell$  and $1 \leq j \leq n-\ell$.
Because the Lie algebra of vector fields preserving a commutative parallelism is commutative (see for instance \cite{kobayashi}), the theorem is proved.
\end{proof}

For Liouville integrable Hamiltonian vector fields, this theorem was already proved along same lines by Ramis \cite{moralesramissimo}.

\section{Linear difference equations}
\label{lineaire}

\subsection{Basic concepts}

For the convenience of the reader, we recall some well-known
concepts from difference equations theory.

\begin{defi}
Let $\phi$ be a non-periodic Moebius transformation of $\mathbb{P}^1(\complex)$
and $E \to \mathbb{P}^1(\complex)$ a rank $n$ vector bundle on the
projective line. A rank $n$ $\phi$-difference system on $E$ is a
lift $\Phi : E \to E$ of $\phi$ linear on the fibers. Let $U$ be a
$\phi$-invariant open set of $\mathbb{P}^1(\complex)$ for the
transcendent topology. A section $Y : U \to E$ is a solution of
$\Phi$ if its graph is $\Phi$-invariant.
 \end{defi}

Using a rational change of variable any Moebius transformation can be
converted either into a translation $\phi : z \mapsto z+h$ with $h
\in \complex$ or into a $q$-dilatation  $\phi : z \mapsto qz$ with
$q\in \complex^*$, depending on the number of fixed points of
$\phi$ on $\mathbb{P}^1(\complex)$ (one in the first case, two in
the second case). The corresponding equations are respectively
called finite difference and $q$-difference equations.\\

By a gauge transformation, any $\phi$-difference system can be
transformed in a system on a trivial vector bundle. In this case,
considering a coordinate $z$ on the projective line, $\Phi$ is
described by an invertible matrix with rational entries :
$$
\begin{array}{lccc}
\Phi : &\mathbb{P}^1(\complex) \times \complex^n &\to& \mathbb{P}^1(\complex) \times \complex^n\\
         & (z,Y) &\mapsto& (\phi z, A(z)Y)
\end{array}
$$
with $A \in GL_{n}(\complex(z))$. The equations satisfy by a
solution are :
\begin{equation}
Y(\phi z) = A(z)Y(z).
\end{equation}

In the sequel, vector bundles will be endowed with a fiberwise symplectic form $\omega$ and $\Phi$ will preserve this form.
By a fiberwise linear change of variables, one can assume that $\omega$ is the canonical symplectic $\J$ form on fibers.
Such a $\phi$-difference system has a matrix $A$ in $Sp_{2n}(\complex(z))$.

What precedes can be rephrased and generalized using the language
of difference modules; these notions being used in \cite{singer-vanderput}, one of
our main references for the Galois theory of difference equations
with \cite{sauloy}, we shall remind them. A \textit{difference ring} is
a couple $(R,\phi)$ where $R$ is a commutative ring and where
$\phi$ is a \textit{difference operator}, that is an automorphism
of $R$. If $R=k$ is a field, then $(k,\phi)$ is called a
\textit{difference field}. A \textit{difference module} $M$ over a
difference field $(k,\phi)$ is a $k[\phi,\phi^{-1}]$-module, whose
$k$-algebra obtained after restriction of scalars is finite
dimensional. The choice of a particular $k$-basis $\mathcal E$ of
$M$ provides a difference system with coefficients in $k$ : $\phi
Y=AY$ where $A\in Gl_n(k)$ is the matrix representing $\phi
\mathcal E$ in $\mathcal E$; the difference system obtained by
choosing another basis $\mathcal E'$ takes the form $\phi Y=A'Y$
with $A'=(\phi P)^{-1}AP$ and $P\in Gl_n(k)$. Conversely, starting
from a difference system $\phi Y=A Y$, $A\in Gl_n(k)$ we construct
a difference module in an obvious way; the difference module
obtained from $\phi Y=A'Y$ with $A'=(\phi P)^{-1}AP$ and $P\in
Gl_n(k)$ is isomorphic to those obtained from $A$.

\subsection{Picard-Vessiot theory for linear difference equations}

In the whole section we consider a difference field $(k,\phi)$
with algebraically closed field of constants (we remind that the
field of constants of $(k,\phi)$ is $C=\{c\in k \ | \ \phi c =c
\}$) and : \begin{equation}\label{syst aux diff}\phi Y=AY
\end{equation} a difference system with coefficients in $k$
(\textit{i.e.} $A\in Gl_n(k)$).

We shall remind the construction of the difference Galois group of (\ref{syst
aux diff}) due to M. van der Put and M. Singer in \cite{singer-vanderput}.

The counterpart for difference equations of the field of
decomposition of an algebraic equation is the Picard-Vessiot
extension.

\begin{defi}
A Picard-Vessiot ring for the difference system (\ref{syst aux
diff}) is a $k$-algebra $R$ such that :
\begin{itemize}
\item[(i)] an automorphism of $R$, also denoted by $\phi$, which extend $\phi$ is given;
\item[(ii)] $R$ is a simple difference ring;
\item[(iii)] there exists a fundamental matrix $\Y$ for $\phi Y=AY$ having entries in $R$ such that $R=k[(\Y_{i,j})_{1 \leq i,j \leq n},\det(\Y)^{-1}]$.
\end{itemize}
\end{defi}

Let us precise the terminology used in the previous definition :
\begin{itemize}
\item[$\bullet$] a difference ideal of a difference ring is an ideal $\I$ stable by $\phi$.
We say that $R$ is a simple difference ring if it has only trivial difference ideal;
\item[$\bullet$] a fundamental matrix $\Y$ for $\phi Y=AY$ having entries in $R$ is a matrix
$\Y \in Gl_n(R)$ such that $\phi \Y=A\Y$.
\end{itemize}

Van der Put and Singer proved in \cite{singer-vanderput} that such an extension
exists and that it is unique up to isomorphism of difference
rings. Furthermore, the field of constants of the Picard-Vessiot
extension is equal to the fields of constants of the base field.

Let us give the general lines of the proof of the existence of the
Picard-Vessiot extension. We denote by $X=(X_{i,j})_{1\leq i,j
\leq n}$ a matrix of indeterminates over $k$ and we extend the
difference operator $\phi$ to the $k$-algebra
$$
U = k\left[(X_{i,j})_{1\leq
i,j \leq n},\det(X)^{-1}\right]
$$
 by setting $(\phi X_{i,j})_{1\leq i,j \leq
n}=A(X_{i,j})_{1\leq i,j \leq n}$. Then, for any maximal
difference ideal $\I$ (which is not necessarily a maximal ideal) of
$U$, the couple
$(R=U/\I,\phi)$ is a
Picard-Vessiot ring for the difference system (\ref{syst aux
diff}).

\begin{defi}
The Galois group $G$ of the difference system (\ref{syst aux
diff}) is the
group of difference automorphisms of its Picard-Vessiot extension
over the base field $(k,\phi)$. It is a linear algebraic group
over the field of constants $C$.
\end{defi}

If $A$ `takes its values' in a
particular linear algebraic group, then the Galois group of
(\ref{syst aux diff}) is a subgroup of this algebraic group. For
later use, let us reprove this statement in the symplectic case.

\begin{lem}\label{PV symp}
Assume that $k=\complex(z)$ and that $A \in Gl_n(\complex(z))$ is
sympletic : $\vphantom{}^t A \J A=Id$ where $\J$ denotes the matrix
of the symplectic form. Then there exists $R$ a Picard-Vessiot
extension of (\ref{syst aux diff}) containing a symplectic
fundamental system of solutions of $\phi Y=AY$, that is there
exists $\Y$ a fundamental system of solutions of $\phi Y=AY$ with
entries in $R$ such that $\vphantom{}^t \Y \J \Y=Id$. The
corresponding Galois group belongs to the symplectic group.
\end{lem}

\begin{proof}As above, we denote by $X=(X_{i,j})_{1\leq i,j
\leq n}$ a matrix of indeterminates over $k$ and we extend the
difference operator $\phi$ to the $k$-algebra
$$U = k[(X_{i,j})_{1\leq i,j \leq n},\det(X)^{-1}]
$$
 by setting $(\phi X_{i,j})_{1\leq i,j \leq
n}=A(X_{i,j})_{1\leq i,j \leq n}$. Remark that the ideal $\tild \I$
of $U$ generated by entries of
$\vphantom{}^t X \J X-Id$ is a difference ideal. Indeed this ideal
is made of the linear combinations $\sum_{i,j}
a_{i,j}(\vphantom{}^t X \J X-Id)_{i,j}$ and
$$
\begin{array}{rcl}
\phi \sum_{i,j}a_{i,j}(\vphantom{}^t X \J X-Id)_{i,j} & = & \sum_{i,j} \phi (a_{i,j})(\vphantom{}^t \phi (X) \J \phi(X)-Id)_{i,j} \\
& = & \sum_{i,j} \phi (a_{i,j}) (\vphantom{}^t (AX) \J (AX)-Id)_{i,j} \\
& = & \sum_{i,j} \phi (a_{i,j})(\vphantom{}^t X \J X-Id)_{i,j}.
\end{array}
$$
 Hence, there
exists a maximal difference ideal $\I$ of $R$ containing $\tild \I$,
so that $U/\I$ is a
Picard-Vessiot ring for $\phi Y=AY$. The lemma follows easily.
\end{proof}

We come back to the general non-necessarily symplectic situation. In contrast with the differential case,
the Picard-Vessiot ring is not necessarily a domain. Let us go further into the study of its
structure. Following \cite{singer-vanderput}, the Picard-Vessiot ring, say $R$,
can be decomposed as a product of domains : $R=R_1 \oplus \cdots
\oplus R_s$ such that :

\begin{itemize}
  \item[(i)] $\phi R_i=R_{i+1}$;
  \item[(ii)] each component $R_i$ is a
Picard-Vessiot ring for the difference equation~:
\begin{equation}\label{equa aux diff bis}\phi^s Y=\left(\phi^{s-1} A \cdots A\right)
Y,\end{equation} the corresponding Galois group is denoted by $G'$.
\end{itemize}

The proof of the following lemma is left to the reader.

\begin{lem}\label{lem integ t}
If the difference system (\ref{syst aux diff}) is integrable then
it is also the case of the difference system (\ref{equa aux diff
bis}).
\end{lem}

Furthermore, following \cite{singer-vanderput}, we have a short
exact sequence :
\begin{equation}\label{exact seq}0\rightarrow G' \rightarrow G \rightarrow \rinteger / s
\rinteger \rightarrow 0.\end{equation}

The following result will be essential in what follows.

\begin{lem}\label{lemme galois inclusion}
If $G'^0$ is abelian, then $G^0$ is abelian.
\end{lem}

\begin{proof}
Indeed, the above exact sequence (\ref{exact seq}) allows us to identify
$G'$ with an algebraic subgroup of $G$ of finite index, so that
$G^0 \subset G'^0$.
\end{proof}

\section{Linearization of integrable dynamical systems}
\label{linearisation}

In this section, starting with a discrete dynamical
system $f$ having an invariant curve of genus zero on which $f$ is
a Moebius transform, we explain how one can obtain by a
linearization procedure a $\phi$-difference system (\textit{finite difference or
$q$-difference system}). Differential invariants of $f$ provide in some cases
differential invariants for its linearization.

\subsection{Adapted curves and discrete variational equations}

\begin{defi}
A curve $\phi$-adapted to $f$ is a rational embedding
$\iota :\mathbb{P}^1(\complex) \dasharrow V$ such that $f \circ
\iota = \iota \circ \phi$.

\end{defi}

\begin{defi}
Let $\iota$ be a curve $\phi$-adapted to $f$. The difference system $\phi Y= Df(\iota) Y$ over
$\mathbb{P}^1(\complex)$ is called the discrete variational
equation of $f$ along $\iota$.
\end{defi}

By a rational gauge transform, the pull back of the tangent bundle : $ \iota^* TV$ can be assume to be trivial and the variational equation can be written as a $\phi$-difference system on the projective line as stated in the definition.

\subsection{Integrability of the discrete variational equation}

In this section we prove that the discrete variational equation of
$f$ along a $\phi$-adapted curve is integrable in the following
sense :

\begin{defi}\label{def integ linear}
A rank $2n$ difference equation $\phi Y=A Y$, $A\in
SL_{2n}(\complex(z))$ is integrable if it admits $n+\ell$ independent rational
first integrals $h_1, \ldots, h_{n +\ell}$ such that the distribution  $\ker dz \cap \bigcap_i \ker dh_i  $
is generated by the fiberwise symplectic gradients $X_{h_1},\ldots X_{h_{n-\ell}}$.
\end{defi}

There are two points requiring comments :
\begin{itemize}
  \item[$\bullet$] a function $H \in \complex(z,Y_{i} | 1\leq i \leq 2n)$ is a
first integral of a difference system $\phi Y=AY$ if $H(\phi z, A
Y)=H(z, Y)$.
\item[$\bullet$] The symplectic structure on the fibers is given by the constant canonical symplectic structure on $\complex^{2n}$.
\end{itemize}

Our main tool is the \textit{generic junior part} that we shall
now introduce.

Let us consider a rational embedding
$\iota:\mathbb{P}^1(\complex) \dasharrow V$. The set of rational
functions on $V$ whose polar locus does not contains
$\iota(\mathbb{P}^1(\complex))$ is denoted by
$\complex[V]_{\iota}$

\begin{defi}
The generic valuation $\val_\iota(H)$ of a function $H \in
\complex[V]_\iota$ along $\iota$ is defined by :
$$\val_\iota(H)=\min \{k\in \integer \ | \ D^kH(\iota) \not
\equiv 0 \text{ over } \mathbb{P}^1(\complex) \}.$$ It extends to
$H\in \complex(V)$ by setting
$\val_\iota(H)=\val_\iota(F)-\val_\iota(G)$, where $H=F/G$
with $F, G \in \complex[V]_\iota$.
\end{defi}

\begin{defi}
 Let us consider $H\in \complex[V]_\iota$ and $T_{\iota} = \iota^* TV$. We define
$H^\circ_\iota : T_{\iota} \dasharrow \complex$, the generic junior part of $H$ along $\iota$,
by :
$$H^\circ_\iota=D^{\val_\iota(H)}F(\iota).$$

We extend this definition to $H\in \complex(V)$ by
setting $H^\circ=F^\circ/G^\circ$ for any $F,G \in \complex[V]_\iota$
such that $H=F/G$.
\end{defi}
The reader will easily verify that the above definitions do not depend on a particular choice of
$F,G$.\\

\begin{rem}
\label{genericandpoint}
\textit{Junior parts of $F$ at points $p \in \iota(\mathbb{P}^1(\complex))$ in the sense
of \cite{morales} coincide with the generic junior part over a
Zariski-dense subset of $\mathbb{P}^1(\complex)$.}
\end{rem}

The interest of the generic junior part in our context consists in
the fact that it converts a first integral of $f$ into a
first integral of the variational equation.

\begin{lem}
Let us consider $H \in \complex(V)$ a first integral
of $f$. Then $H^\circ_\iota$ is a first integral of the
variational equation of $f$ along $\iota$.
\end{lem}

\begin{proof}
Let us consider $F,G \in \complex[V]_\iota$ such that $H=F/G$. The
assertion follows by differentiating $v_\iota(F)+v_\iota(G)$
times the equality $(F\circ f) \cdot G=(G\circ f) \cdot F$.
\end{proof}

Therefore, if $H_1,...,H_{n+\ell}$ are $n+\ell$ first integrals of
$f$, then
$(H_{1})_{\iota}^\circ,\ldots,(H_{n+\ell})^\circ_{\iota}$ are
$n+\ell$ first integrals of the variational equation of $f$ along
$\iota$. However, in general, these function are \textit{not}
functionally independent, even if $H_1,...,H_{n+\ell}$ are
functionally independent. This difficulty can be overcome by using
the following lemma due to Ziglin (\cite{ziglin}).

\begin{lem}
Let $F_1,\ldots,F_k$ be $k$ meromophic functions functionally
independent in a neighborhood of $\iota(\mathbb{P}^1(\complex))$. Then there exist $k$
polynomials $P_1,\ldots,P_{k} \in \complex[z_1,\ldots,z_{k}]$
such that
the generic junior part $(G_{1})_{\iota}^{\circ}
,\ldots,(G_{k})_{\iota}^{\circ}$ along $\iota$ of $G_{1} =
P_1(F_1,\ldots,F_k),\ldots,G_{k} = P_k(F_1,\ldots,F_k)$ are
functionally independent.
\end{lem}

\begin{proof}
The usual Ziglin Lemma ensures that the assertion is true if one
replaces `generic junior part'  by `junior part at a given point'. The previous remark allows us to conclude the proof.
\end{proof}

\begin{theo}\label{integ de la varia}
The variational equation of $f$ along a $\phi$-adapted curve
$\varphi$ is integrable.
\end{theo}
\begin{proof}Let $H_{1},\ldots,H_{n+\ell}$ be first integrals of $f$ with symplectic gradients $X_{H_{1}},\ldots,X_{H_{n-\ell}}$
generating the distribution $\bigcap_{i} \ker dH_{i}$.

By Ziglin lemma, junior parts of the first $n-\ell$ first
integrals give $n-\ell$ independent first integrals
$(H_{1})_{\iota}^{\circ},\ldots,(H_{n-\ell})_{\iota}^{\circ}$ for
the variational equation.  The independence of the differentials
$d(H_{1})_{\iota}^{\circ},\ldots,d(H_{n-\ell})_{\iota}^{\circ}$
implies the independence of their symplectic gradients
$X_{(H_{1})_{\iota}^{\circ}},\ldots,X_{(H_{n-\ell})_{\iota}^{\circ}}$.
The set of first integrals is completed applying Ziglin lemma to
the $2 \ell$ remaining first integrals.
\end{proof}

In the next sections we will prove that the integrability of the
variational equation implies severe constraints on its Galois
group (more precisely on the neutral component of the Galois
group).

\section{Integrability and Galois theory}
\label{theoreme}

\subsection{Statement of the main theorem}

\begin{theo} \label{main theo}
Let $f$ be an integrable symplectic map having a
$\phi$-adapted curve $\iota$. Then the neutral component of the
difference Galois group of the discrete variational equation of $f$ along
$\iota$ is commutative.
\end{theo}

In the following sections we give two proofs of our main theorem :
the first one is based on Picard-Vessiot theory developed by Van
der Put and Singer, the second one is based on Malgrange groupoid.

\subsection{First proof of the main theorem : Picard-Vessiot approach}

\begin{lem}\label{tild}
Let $(k=\complex(z),\phi)$ be a difference field with field of
constants $\complex$. Let $H$ be a first integral of the
difference system with coefficients in $\complex(z)$ : $\phi Y=A
Y$. Let $R$ be a Picard-Vessiot extension for this difference
system which is supposed to be a domain. Let $\Y \in Gl_n(R)$ be a
fundamental system of solutions with coefficients in $R$. Then the
function~:
\begin{eqnarray*}
\vphantom{}_\Y H : \complex^n& \rightarrow & \complex\\
y& \mapsto & H(z,\Y y)
\end{eqnarray*}
is invariant under the action of the Galois group $G$.
\end{lem}

The action of the Galois group is given by the following : for
$(\sigma,y) \in G\times \complex^n$, we set $\sigma \cdot
y=C(\sigma) y$ where $C(\sigma)\in Gl_n(\complex)$ is such that
$\sigma \Y=\Y C(\sigma)$.

\begin{proof} The first point to check is that, for any
$y\in \complex^n$, $H(z,\Y y)\in\complex$. \textit{A priori} $H(z,\Y y)$
belongs to $\text{Frac}(R)$, the Picard-Vessiot field, but $\phi
H(z,\Y y)= H(\phi z, (\phi\Y) y) = H(\phi z, (A(z) \Y) y)= H(\phi
z, A(z) (\Y y))=H(z,\Y y)$ so that $H(z,\Y y)$ is $\phi$-invariant
: $H(z,\Y y)$ belongs to the field of contant of the
Picard-Vessiot field, that is to $\complex$. Since $H(z,\Y y)$
lies in the field of constants $\complex$, it is invariant by the
Galois group : $\sigma H(z,\Y y)=H(z,\Y y)$. But, since $H$ has
rational coefficients in $z$ and that $\complex(z)$ is invariant
by the action of the Galois group, we also have, for all $\sigma
\in G$, $\sigma H(z,\Y y)=H(z,\sigma \Y y)=H(z,\Y
C(\sigma)y)=H(z,\Y (\sigma \cdot y))$.
\end{proof}

\begin{lem}\label{lem integre}
Suppose that the symplectic difference equation $\phi Y=AY$ with
coefficients in $\complex(z)$ is integrable and that its
Picard-Vessiot ring is a domain. Then the Lie algebra of its
Galois group is abelian.
\end{lem}

\begin{proof}
Let $H_1,...,H_{n+ \ell}$ be $n+ \ell$ first integrals insuring integrability of the difference equation as in definition \ref{def integ linear}. Lemma \ref{PV symp} ensures that there exists a
Picard-Vessiot extension $R$ containing a symplectic fundamental
system of solution $\Y \in \text{Sp}_{2n}(R)$.
Because
$$
\sum_{j} \frac{\partial \vphantom{}_\Y H_i}{\partial y_{j}} dy_{j} = \sum_{j,k} \frac{\partial  H_i}{\partial Y_{k}} \Y_{k,j} dy_{j}
$$
and invertibility of $\Y$ independence of the $n+\ell$ forms $\sum_{k} \frac{\partial H_{i}}{\partial Y_{k}}dY_{k}$ implies independence the $n+ \ell$ differentials $d\vphantom{}_\Y H_i$.
Moreover from lemma \ref{tild} this functions are invariant under the action of the Galois group.

Let $X$ be a linear vector field on $\complex^{2n}$ in the Lie algebra of the Galois group.

By the assumption of integrability, all vector field preserving $\vphantom{}_\Y H_1,...,\vphantom{}_\Y H_{n+\ell}$ are combinations of the $n-\ell$ first symplectic gradients:

$$
X = \sum_{i=1}^{n-\ell} c_{i} X_{\vphantom{}_\Y H_i}
$$
with $c_{i} \in \complex(y_{1}, \ldots, y_{2n})$.

Because the Galois group is a subgroup of $\text{Sp}_{2n}(\complex)$, if it preserves a function, it  preserves its symplectic gradient too. Thus one gets  $[X_{\vphantom{}_\Y H_j}, X] = 0$ for all $1\leq j \leq n + \ell$ and
$$
dc_{i}  = \sum_{i=1}^{n + \ell} c_{i}^j d\vphantom{}_\Y H_j
$$
for some $c_{i}^j \in \complex(y_{1}, \ldots, y_{2n})$.
 If $Y$ is a second vector field in the Lie algebra of the Galois group, an easy computation shows that $[X,Y] = 0$. This proves the lemma
\end{proof}

\begin{prop}\label{prop ab gen}
If a linear symplectic difference system is integrable then the
neutral component of its Galois group is abelian.
\end{prop}

\begin{proof} Invoking lemma \ref{lem integ t}, we get that for every
integer $s$ the system $\phi Y=\phi^{s-1}A \cdots \phi A \cdot A Y$
is symplectic and integrable. The lemma follows from this
observation and from lemma \ref{lemme galois inclusion} together
with lemma \ref{lem integre} (indeed maintaining the notations of
lemma \ref{lemme galois inclusion}, lemma \ref{lem integre}
implies that $G'^0$ is commutative because $G'$ is the Galois
group of a symplectic and integrable system having a domain
($R_i$) as a Picard-Vessiot ring).
\end{proof}

\noindent \textit{First proof of theorem \ref{main theo}.\quad} It is a direct consequence of Theorem \ref{integ de la varia} and
Proposition \ref{prop ab gen}. \hfill $\square$

\subsection{Second proof using Malgrange groupoid}

Let $f : V \dasharrow V$ be a rational map  and $\iota : \mathbb{P}^1(\complex) \dasharrow V$ a $\phi$-adapted curve.
\begin{lem}
\label{21} If the Malgrange groupoid of $f$ is infinitesimally
commutative then the Malgrange groupoid of $R_{q} f$ is
infinitesimally commutative too.
\end{lem}
\begin{proof}
Differential invariants of $f$ give differential invariants for $R_{q}f$  as described in \ref{invariantsdifferentiels}. Because $R_{q}f$ is a prolongation it preserves also the vector field $D_{j}$ on $R_{q}V$ : $(R_{q}f)^*D_{j} = D_{j}$. Hence Malgrange groupoid of $R_{q}f$ must preserve these vector fields. A direct computation shows that a vector field commutes with all the $D_{j}$ if and only if it is a prolongation of a vector field on $V$. The compatibility of the prolongation with the Lie bracket proves the lemma.
\end{proof}

Let $B_{q}\iota$ be the pull-back by $\iota$ of the frames bundle $R_{q}V$. The principal $GL_{2n}(\complex)$-bundle  $B_{1}\iota$ is the bundle of frames of $\iota^*TV$ on $\mathbb{P}^1(\complex)$. Let $R_{q}\iota$ denote the extension of $\iota$ from $B_{q}\iota$ to $R_{q}V$. The system given by the restriction of $R_{1}f$ to $B_{1}\iota$ is the fondamental discrete variational equation. Its solutions are the fondamental solutions of the variational equation. The following lemma follows from the definition of the discrete variational equation.

\begin{lem}
\label{22}
The pull-back by $R_{1}\iota$ of differential invariants defined in a neighborhood of $\iota(\mathbb{P}^1(\complex))$ are differential invariants of the fundamental variational equation.
\end{lem}
We will see in the second proof of  theorem \ref{main theo} that some differential invariants are less easy to grab.
\begin{lem}
\label{malglineaire}
If the difference system $\phi Y = A Y$ on a trivial $SL_{2n}(\complex)$-principal bundle on $\mathbb{P}^1(\complex)$ with $A \in SL_{2n}(\complex(z))$ is integrable then the fibers-tangent vector fields of the Lie algebra of its Malgrange groupoid commute.
\end{lem}
\begin{proof}
The proof is follows exactly the proof of theorem \ref{malgcommutatif}
\end{proof}

Our second proof of theorem \ref{main theo} rely on the following conjectural result :\\

\noindent \textbf{Conjectural Statement.} \textit{For linear ($q$-)difference systems, the action of Malgrange groupoid on the fibers gives the classical Galois groups.}\\

The similar statement for linear differential equations is \textit{true} and is proved in \cite{malgrange}. The proof of this conjectural statement in the $q$-difference case is a work in progress by A. Granier.\\

\noindent \textit{Second proof of theorem \ref{main theo}.\quad}
Lemmas \ref{21} and \ref{22} prove the theorem when the first
integrals are independent on a neighborhood of the invariant
curve.

Let us first assume that the first integrals are defined on
$\iota$ {it i.e.} are in $\complex[V]_{\iota}$ but not independent
on the curve. Let $H$ be a first integral defined in a
neighborhood of the invariant curve and $k$ be its generic
valuation. The form $dH$ is invariant and gives $n$ differential
invariants of order $1$. If $k>1$ then these invariant vanish
above the invariant curve. The differentials of these invariants
induced order $2$ differential invariants. After $k$ iteration of
this process one gets order $k$ differential invariants but theirs
restrictions on $B_{k}\iota$ give well-defined functions on
$B_{1}\iota$. These functions are induced on $B_{1}\iota$ from the
junior part $H_{\iota}^\circ$ on $\iota^*TV$.

By Ziglin lemma one can assume that the junior parts of the first integrals hence the order $k$ invariant induced are independents. One can apply Lemma \ref{malglineaire} to get the infinitesimal commutativity of the vertical part of the Malgrange groupoid and by the above conjectural statement the almost commutativity of the difference Galois group.\\

Now let the first integral be any kind of rational functions on $V$. Such a $H$ can be written $H = \frac{F}{G}$ with $F$ and $G$ in $\complex[V]_{\iota}$. These two functions induced functions on the frame bundles which are no more invariant. As in the previous case, the junior parts $F^\circ$ and $G^\circ$ of $F$ and $G$ are the restrictions of functions induced by $dF$ and $dG$ on the frame bundle. Let denote them by $h^i_{F}$ dans $h^i_{G}$ for $1\leq i \leq n$. Because $H$ is invariant, $F$ and $G$ are semi invariants : $F\circ f = a F$ and $G \circ f = a G$. Let's have a look at the lift of these functions on $R_{1}V$. One has $(h_{F}^{j}) (R_{1}f) = h_{da}^{j} F + a (h_{F}^{j})$. Let $k$ be the generic valuation of $F$ along $\iota$ then after $k$ lift of the differential as functions on higher frame bundles one gets functions $h^{\alpha}_{F}$ such that $h^\alpha_{F} (R_{k}f) = \ldots + a h^\alpha_{F}$ where the dots stand for terms vanishing above $\iota(\mathbb{P}^1(\complex))$. The same computation on $G$ shows that the junior part of $H$ is an invariant of $R_{1}f$ above $\iota$. The arguments used to conclude in the first case can be apply in this situation. This proves the theorem.
\hfill $\square$

\section{Some examples}

In this section we present some examples in the case $2n = 2$. In this situation rational maps $f : V \dasharrow V$ are allowed to be non-symplectic and integrability of $f$ means existence of a non-constant rational first integral.

\subsection{Variants when $n=1$}
Following same lines as the proof of the main theorem, one gets.

\begin{theo}
Let $V$ be an algebraic variety of dimension $2$ and consider a rational map $f : V \dasharrow V$ having a
$\phi$-adapted curve $\iota$. If $f$ gets a rational first integral then the intersection of the neutral component of the
Galois group of the discrete variational equation of $f$ along
$\iota$ with $\text{SL}_2(\complex)$ is commutative.
\end{theo}

\begin{proof}
The proof follows the same lines as that of the main theorem (mores simple actually).
\end{proof}

We also have the following useful result.

\begin{theo}
Let $V$ be an algebraic variety of dimension $2$ and consider a rational function$f : V \dasharrow V$having a
$\phi$-adapted curve $\iota$. If one of the following conditions relative to the Galois group $G$ of the discrete variational equation of $f$ along
$\iota$ holds then $f$ is not integrable :
\begin{itemize}
\item[-] $\text{diag}(\complex^*,\complex^*)$ is conjugated to a subgroup of $G$;
\item[-] $\pmatrice{1}{\complex}{0}{\complex^*}$ is conjugated to a subgroup of $G$.
\end{itemize}
\end{theo}

\begin{proof}
Use the same arguments as for the proof of the main theorem and the fact that a rational function $F(z)(x,y)$ invariant by the natural action of the above groups (on $(x,y)$) is necessarily constant.
\end{proof}

\subsection{Example 1}
Consider
$$
\begin{array}{rccc}
f : & \complex^2 & \dasharrow & \complex^2\\
 & \begin{pmatrix} x \\ y \end{pmatrix} &\mapsto &\begin{pmatrix}qx \frac{1}{1-\frac{y}{x-1}} \\ \underline{q}y
(1-\frac{y}{x-1}))\end{pmatrix}
\end{array}
$$
where $q,\underline{q}\in \complex^*$. The
curve $\varphi(x)=(x,0)$ is a $\sigma_q$-adapted to $f$. The
variational equation of $f$ along $\varphi$ is given by :
\begin{equation}\label{ex5}Y(qx)=\pmatrice{q}{q \frac{x}{x-1}}{0}{\underline{q}}Y(x).\end{equation}
If $\underline{q}$ is in general position then its Galois group is
$\pmatrice{1}{\complex}{0}{\complex^*}$ so that $f$ is not
integrable. \textit{A contrario}, remark that when
$q\underline{q}$ is a root of the unity then the map $f$ is
integrable ($F(x,y)=(xy)^n$ is a first integral for $n \in
\integer$ large enough).

\subsection{Example 2}

Consider
$$
\begin{array}{rccc}
f : & \complex^2 & \dasharrow & \complex^2\\
 & \begin{pmatrix} x \\ y \end{pmatrix} &\mapsto &\begin{pmatrix}qx + y a(x) \\ yb(x)\end{pmatrix}
\end{array}
$$
with $b\neq 0$. The curve $\varphi(x)=(x,0)$ is a
$\sigma_q$-adapted to $f$. The variational equation of $f$ along
$\varphi$ is given by :
\begin{equation}\label{ex3}Y(qx)=\pmatrice{q}{a(x)}{0}{b(x)}Y(x).\end{equation} Let
us for instance study the case where
$a(x)=b(x)=q\frac{x-1}{bx-c/q}$ with $b/c\not \in q^\rinteger$.
Remark that
$$
\pmatrice{q}{q\frac{x-1}{bx-c/q}}{0}{q\frac{x-1}{bx-c/q}}=
qz\pmatrice{0}{1}{-1}{1}\pmatrice{0}{1}{-\frac{z-1}{bz-c/q}}{\frac{(1+b)z-(1+c/q)}{bz-c/q}}z^{-1}\pmatrice{0}{1}{-1}{1}^{-1},
$$
so that (\ref{ex3}) is rationally equivalent to (the system
associated to) a basic hypergeometric equation with parameters
$(1,b,c)$, the Galois group of which is
$\pmatrice{1}{\complex}{0}{\complex^*}$ (see \cite{roques galois
basic hyper}). Consequently, $f$ in not integrable.

\subsection{Example 3}

Consider
$$
\begin{array}{rccc}
f : & \complex^2 & \dasharrow & \complex^2\\
 & \begin{pmatrix} x \\ y \end{pmatrix} &\mapsto &\begin{pmatrix} ya(x)-x+1 \\ xyb(y)\end{pmatrix}
\end{array}
$$
where $a\in\complex(x)$ and $b \in \complex(y)$ without pole in
$0$. The curve $\varphi(x)=(x,0)$ is a $\tau_{-1}$-adapted curve
with $\tau_{-1}(z)=z-1$. The variational equation of $f$ along
$\varphi$ is given by :
$$Y(1-x)=\pmatrice{-1}{a(x)}{0}{xb(0)}Y(x).$$
The algorithm exposed in \cite{hendricks} allows us to exhibit
many non integrable cases.

\subsection{Example 4}
Consider
$$
\begin{array}{rccc}
f : & \complex^2 & \dasharrow & \complex^2\\
 & \begin{pmatrix} x \\ y \end{pmatrix} &\mapsto &\begin{pmatrix}qx + \tilde{q}y + (x-y) a(x,y) \\ \underline{q}x + \underline{\tilde{q}}y + (x-y) \underline
{a}(x,y))\end{pmatrix}
\end{array}
$$
where $a,\underline{a}\in\complex(x)$ and $q,\tilde{q},\underline{q},\underline{\tilde{q}} \in \complex$
 are such that $\mathfrak{q}=q+\tilde{q}=\underline{q}+\underline{\tilde{q}}\in \complex^*$.
The curve $\varphi(x)=(x,x)$ is a
$\sigma_\mathfrak{q}$-adapted to $f$. The variational equation of $f$ along
$\varphi$ is given by :
\begin{equation}\label{ex4}Y(\mathfrak{q}x)=\pmatrice{q+a(x,x)}{\tilde{q}-a(x,x)}{\underline{q}+\underline
{a}(x,x)}{\underline{\tilde{q}}-\underline
{a}(x,x)}Y(x).\end{equation}

This system is rationally equivalent to :
\begin{equation}\label{ex4}Y(\mathfrak{q}x)=\pmatrice{q-\underline{q}+a(x,x)-\underline{a}(x,x)}{0}{\underline{q}+\underline{a}(x,x)}{\mathfrak{q}}Y(x).\end{equation}

One can see that generally $f$ is not integrable.

\end{document}